\input amstex
\documentstyle{amsppt}
%
%
\nopagenumbers
\pagewidth{360pt}
\pageheight{606pt}
\def\negskp{\hskip -2pt}
\def\sqbi#1{{\sssize[\ssize\hskip 0.5pt #1\sssize]}}
\rightheadtext{Orthogonal matrices in composing tests.}
\topmatter
\title
Orthogonal matrices with rational components
in composing tests for High School students.
\endtitle
\author
Ruslan~A\.~Sharipov
\endauthor
\address Rabochaya~str\.~5, 450003, Ufa, Russia
\endaddress
\email \vtop to 20pt{\hsize=280pt\noindent
R\_\hskip 1pt Sharipov\@ic.bashedu.ru\newline
ruslan-sharipov\@usa.net\vss}
\endemail
\urladdr
http:/\negskp/www.geocities.com/CapeCanaveral/Lab/5341
\endurladdr
\abstract
Fermat Last Theorem, which inspired mathematicians during 300 years,
is proved by Andrew Wiles. Even among mathematicians there is a narrow
circle of specialists, who can read this proof and understand all details.
Is it a reason for pessimism\,? No, since arithmetics if entire numbers
contains broad variety of problems with a simple statement, which might be
not less intricate. One of them arises in elementary geometry.
\endabstract
\endtopmatter
\document
\head
1. Elementary problem on pyramid.
\endhead
    Primary education (Elementary School) and secondary education
(High School) in Russia are united into one stage that now lasts
11 years (from 6 year old to 17 years old). Mathematics is among
disciplines studied during these years. Below we consider a problem,
which can be suggested to 10-th or 11-th year students in the course
of geometry. It is typical, though is a little more complicated
than usual.
\proclaim{\quad Problem on pyramid}\parshape 1 180pt 180pt In triangular
pyramid $ABCD$ three sides \vadjust{\vskip -15pt\hbox to 0pt{\kern 1pt
\hbox{\special{em:graph 943-1a.gif}}\hss}\vskip 15pt} of triangle $ABC$
in its base are given:
$$
\xalignat 3
&|BC|=a,&&|CA|=b,&&|AB|=c.
\endxalignat
$$
From corners $A$ and $B$ two perpendiculars are drown to the faces $BCD$
and $ACD$ respectively. Their lengths are given:
$$
\xalignat 2
&|AF|=f,&&|BG|=g.
\endxalignat
$$
Find the length of the segment $[FG]$ connecting feet of these two
perpendiculars $[AF]$ and $[BG]$.
\endproclaim
\parshape 4 180pt 180pt 180pt 180pt 180pt 180pt 0pt 360pt
     Let's consider in brief the steps leading to the solution of this
problem. First we draw
all three heights in triangle $ABC$: these are segments $[AH]$, $[BK]$,
$[CM]$. It's known that they cross at one point. Denote it by $L$.
For the sake of simplicity we consider the case when triangle $ABC$ is
acute-angled. In this case point $L$ lies inside the triangle $ABC$.
Now let's apply Pythagor's theorem to rectangular triangles
$AMC$ and $BMC$. As a result we obtain the system of equations with
respect to the length of the segment $[AM]$:
$$
\cases |AM|^2+|CM|^2=|AC|^2,\\
(|AB|-|AM|)^2+|CM|^2=|BC|^2.
\endcases
\hskip -2em
\tag1.1
$$
Solving the system of equations \thetag{1.1}, for lengths of $[AM]$ and
$[BM]$ we get
$$
\gathered
|AM|=\frac{|AB|^2+|AC|^2-|BC|^2}{2\,|AB|},\\
\vspace{1ex}
|BM|=\frac{|AB|^2+|BC|^2-|AC|^2}{2\,|AB|}.
\endgathered
\hskip -2em
\tag1.2
$$
Similar formulas can be obtained for $|AK|$, $|KC|$, $|BH|$, and $|HC|$:
$$
\gather
\gathered
|AK|=\frac{|AC|^2+|AB|^2-|BC|^2}{2\,|AC|},\\
\vspace{1ex}
|CK|=\frac{|AC|^2+|BC|^2-|AB|^2}{2\,|AC|},
\endgathered\hskip -2em
\\
\vspace{1ex}
\gathered
|BH|=\frac{|BC|^2+|AB|^2-|AC|^2}{2\,|BC|},\\
\vspace{1ex}
|CH|=\frac{|BC|^2+|AC|^2-|AB|^2}{2\,|BC|}.
\endgathered\hskip -2em
\endgather
$$
Let's replace $|AB|-|AM|$ by $|BM|$ in the second equation of the system
\thetag{1.1}. Then we can derive the following formula for the length
of segment $[CM]$:
$$
|CM|=\sqrt{\frac{|AC|^2+|BC|^2-|AM|^2-|BM|^2}{2}}.\hskip -2em
\tag1.3
$$
Similar formulas can be derived for the lengths of segments $[AH]$ and
$[BK]$:
$$
\gathered
|AH|=\sqrt{\frac{|AB|^2+|AC|^2-|CH|^2-|BH|^2}{2}},\\
\vspace{1ex}
|BK|=\sqrt{\frac{|AB|^2+|BC|^2-|AK|^2-|CK|^2}{2}}.
\endgathered\hskip -2em
\tag1.4
$$
In order to calculate lengths of segments $[KL]$ and $[HL]$ we use
similarity of triangles: $\triangle KLC\sim\triangle MAC$ and
$\triangle HLC\sim\triangle MBC$. This yields:
$$
\pagebreak
\xalignat 2
&|KL|=\frac{|AM|}{|CM|}\,|KC|,
&&|HL|=\frac{|BM|}{|CM|}\,|HC|.
\hskip -2em
\tag1.5
\endxalignat
$$\par
     Now let's draw segments $[FH]$ and $[GK]$. According to the
theorem on three perpendiculars, we have $GK\perp AC$ and $FH
\perp BC$. Then, since we already know $|AF|$ and $|BG|$, we can
calculate lengths of segments $[FH]$ and $[GK]$:
$$
\gathered
|FH|=\sqrt{|AC|^2-|CH|^2-|AF|^2},\\
\vspace{1ex}
|GK|=\sqrt{|BC|^2-|CK|^2-|BG|^2}.
\endgathered
$$
In order to derive first of these two expressions we applied Pythagor's
theorem to rectangular triangles $AHC$ and $AHF$. Second expression is
derived by Pythagor's theorem applied to triangles $BKC$ and $BKG$.
\par
     Orthogonal projections of the points $F$ and $G$ onto the plane of
the base of pyramid belong to the straight lines $AH$ and $BK$. Denote
these projections by $\tilde F$ and $\tilde G$ respectively. For the sake
of simplicity we consider the case when points $F$ and $G$ are above
the base of pyramid (i\.~e\. in upper halfspace separated by the plane
$ABC$), and when their projections $\tilde F$ and $\tilde G$ belong to
the segments $[HL]$ and $[KL]$ respectively (see Fig\.~1.2 and Fig\. 1.3).
\vadjust{\vskip 15pt\hbox to 0pt{\kern 12pt\hbox{\special{em:graph
943-1b.gif}}\hss}\vskip 132pt}Due to similarity of triangles $\triangle
GK\tilde G\sim\triangle KBG$ and $\triangle HF\tilde F\sim\triangle FAH$
we derive the following formulas:
$$
\xalignat 2
&|G\tilde G|=\frac{|BG|\,|GK|}{|BK|},
&&|F\tilde F|=\frac{|AF|\,|FH|}{|AH|},\hskip -2em
\tag1.6\\
\vspace{1ex}
&|K\tilde G|=\frac{|GK|^2}{|BK|},
&&|H\tilde F|=\frac{|FH|^2}{|AH|}.\hskip -2em
\tag1.7
\endxalignat
$$
The length of the segment $[\tilde G\tilde F]$ (see Fig\.~1.4 below) is
determined by cosine theorem applied to the triangle $\tilde GL
\tilde F$:
$$
|\tilde G\tilde F|=\sqrt{|L\tilde G|^2+|L\tilde F|^2-2\,|L\tilde G|\,
\,|L\tilde F|\,\cos(\widehat{KLH})}\hskip -2em
\tag1.8
$$
Note that angles $\angle KLH$ and $\angle KCH$ complete each other to
a straight angle. Indeed, triangle $LKC$ is rectangular (see Fig\.~1.4
below). Same is true for triangle $LHC$. Hence for the angles of these
two triangles we can write the equalities:
$$
\aligned
&\widehat{KLC}+\widehat{KCL}=90^\circ,\\
&\widehat{HLC}+\widehat{HCL}=90^\circ.
\endaligned
$$
Adding these equalities and taking into account that $\widehat{KLC}
+\widehat{HLC}=\widehat{KLH}$ and $\widehat{KCL}+\widehat{HCL}
=\widehat{KCH}$, \vadjust{\vskip 15pt\hbox to 0pt{\kern 3pt\hbox{
\special{em:graph 943-1c.gif}}\hss}\vskip 200pt}we get the required
equality $\widehat{KLH}+\widehat{KCH}=180^\circ$. As an immediate
consequence of this equality we can write the equality for cosines:
$$
\cos(\widehat{KLH})=-\cos(\widehat{KCH}).
$$
Cosine of the angle $\widehat{KCH}$ can be determined by applying
cosine theorem to the triangle $ABC$, which lies in the base of
pyramid $ABCD$:
$$
\cos(\widehat{KCH})=\frac{|AC|^2+|BC|^2-|AB|^2}{2\,|AC|\,|BC|}.
$$
Lengths of segments $[L\tilde F]$ and $[L\tilde G]$ in formula
\thetag{1.8} can be calculated as follows:
$$
\xalignat 2
&|L\tilde F|=|LH|-|H\tilde F|,
&&|L\tilde G|=|LK|-|K\tilde G|.
\endxalignat
$$
This is obvious from Fig\.~1.2 and Fig\.~1.3. Now the length of
segment $[FG]$, which was to be found, is calculated by Pythagor's
theorem (see Fig\.~1.5):
$$
|FG|=\sqrt{|\tilde F\tilde G|^2+(|G\tilde G|-|F\tilde F|)^2}.
$$\par
     So {\bf problem on pyramid} is solved. This is typical stereometric
problem that can be used to test the knowledge of some basic facts and
spatial imagination of students. Its solution considered just above is not
tricky. But it is rather huge, and we cannot write
simple explicit formula expressing $|FG|$ through parameters $a$, $b$,
$c$, $f$, and $g$. Therefore we should give numeric values for these
parameters, choosing them so that they provide simple numeric values
for ultimate result and for results of all intermediate calculations.
Thus, another problem arises, problem of choosing proper numeric values
for $a$, $b$, $c$, $f$, and $g$. We shall consider this problem below.
\head
2. Orthogonal matrices.
\endhead
     Let's apply the coordinate method to the problem on pyramid. Here
we have two natural triples of orthogonal vectors. First consists of
vectors $\overrightarrow{AF}$, $\overrightarrow{FH}$, $\overrightarrow{HC}$,
second is formed by vectors $\overrightarrow{BG}$, $\overrightarrow{GK}$,
and $\overrightarrow{KC}$. Let's consider three unitary vectors $\bold e_1$,
$\bold e_2$, $\bold e_3$ directed along vectors $\overrightarrow{AF}$,
$\overrightarrow{FH}$, and $\overrightarrow{HC}$. Then choose other three
unitary vectors directed along vectors $\overrightarrow{BG}$,
$\overrightarrow{GK}$, and $\overrightarrow{KC}$. Vectors $\bold e_1$,
$\bold e_2$, $\bold e_3$ and $\bold h_1$, $\bold h_2$, $\bold h_3$ form
two bases consisting of unitary vectors orthogonal to each other. Such
bases are called {\bf orthonormal bases} (ONB). Let's consider the
following expansions binding vectors of two ONB's:
$$
\xalignat 2
&\bold h_i=\sum^3_{i=1}S^k_i\,\bold e_k,
&&\bold e_k=\sum^3_{i=1}T^j_k\,\bold h_j.
\hskip -2em
\tag2.1
\endxalignat
$$
Coefficients of the expansions \thetag{2.1} are usually arranged into
square matrices, which are called {\bf transition matrices}:
$$
\xalignat 2
&S=\Vmatrix
S^1_1 & S^1_2 & S^1_3\\ \vspace{1ex}
S^2_1 & S^2_2 & S^2_3\\ \vspace{1ex}
S^3_1 & S^3_2 & S^3_3
\endVmatrix,
&&T=\Vmatrix
T^1_1 & T^1_2 & T^1_3\\ \vspace{1ex}
T^2_1 & T^2_2 & T^2_3\\ \vspace{1ex}
T^3_1 & T^3_2 & T^3_3
\endVmatrix.\hskip -2em
\tag2.2
\endxalignat
$$
Matrices $S$ and $T$ in \thetag{2.2} implement direct and inverse
transitions from base to base, they are inverse to each other, i\.~e\.
their product is a unitary matrix:
$$
S\cdot T=T\cdot S=E.
$$
If we treat \thetag{2.1} as transition from the base $\bold e_1$,
$\bold e_2$, $\bold e_3$ to the base $\bold h_1$, $\bold h_2$,
$\bold h_3$, then $S$ is called {\bf direct transition matrix},
while $T$ is called {\bf inverse transition matrix}.\par
    Matrices $S$ and $T$ in our case are binding two ONB's. Therefore
components of these matrices are bound by a series of relationships.
If we denote by $S^t$ and $T^t$ transposed matrices, i\.~e\. if we
denote
$$
\xalignat 2
&S^t=\Vmatrix
S^1_1 & S^2_1 & S^3_1\\ \vspace{1ex}
S^1_2 & S^2_2 & S^3_2\\ \vspace{1ex}
S^1_3 & S^2_3 & S^3_3
\endVmatrix,
&&T^t=\Vmatrix
T^1_1 & T^2_1 & T^3_1\\ \vspace{1ex}
T^1_2 & T^2_2 & T^3_2\\ \vspace{1ex}
T^1_3 & T^2_3 & T^3_3
\endVmatrix,
\endxalignat
$$
then these relationships for components of $S$ and $T$ can be written as
follows:
$$
\xalignat 2
&S^t\cdot S=E,&&T^t\cdot T=E.\hskip -2em
\tag2.3
\endxalignat
$$
From \thetag{2.3} and from $S\cdot T=T\cdot S=E$ we immediately derive
$S^t=T$ and $T^t=S$.\par
     Matrices that satisfy the relationships \thetag{2.3} are called
{\bf orthogonal matrices}. Sum of squares of elements in each column
and in each string of orthogonal matrix is equal to $1$. So we have
the relationships
$$
\sum^3_{i=1}(S^i_k)^2=\sum^3_{i=1}(S^k_i)^2=1
\text{\ \ for all \ }k=1,\,2,\,3.\hskip -2em
\tag2.4
$$
Sums of products of elements from different columns and/or different
string are equal to zero. This property is expressed by the relationships
$$
\sum^3_{i=1}S^i_k\,S^i_q=\sum^3_{i=1}S^k_i\,S^q_i=0
\text{\ \ for \ }k\neq q.\hskip -2em
\tag2.5
$$
The relationships \thetag{2.4} and \thetag{2.5} are easily derived from
\thetag{2.3}. Moreover, from \thetag{2.3} one can derive the following
relationships for determinants of $S$ and $T$:
$$
\xalignat 2
&(\det S)^2=1,&&(\det T)^2=1.
\endxalignat
$$
Therefore $\det S=\det T=\pm\,1$. Looking attentively at Fig\.~1, one
can note that $\overrightarrow{AF}$, $\overrightarrow{FH}$,
$\overrightarrow{HC}$ and $\overrightarrow{BG}$, $\overrightarrow{GK}$,
$\overrightarrow{KC}$ are oppositely oriented triples of vectors: first
is {\bf left}, while second is {\bf right}. Hence bases $\bold e_1$,
$\bold e_2$, $\bold e_3$ and $\bold h_1$, $\bold h_2$, $\bold h_3$ are
also oppositely oriented. This fact is reflected by the sign of
determinants of transition matrices:
$$
\det S=\det T=-1.\hskip -2em
\tag2.6
$$\par
    Further we shall be interested in the case when all components of
matrices $S$ and $T$ are rational numbers. Components of $S$ belonging
to the same column can be brought to common denominator, and hence, they
can be written as
$$
\xalignat 3
&S^1_1=\frac{p_1}{d_1},&&S^2_1=\frac{p_2}{d_1},&&S^3_1=\frac{p_3}{d_1},
\hskip -2em\\
\vspace{1ex}
&S^1_2=\frac{q_1}{d_2},&&S^2_2=\frac{q_2}{d_2},&&S^3_2=\frac{q_3}{d_2},
\hskip -2em\tag2.7\\
\vspace{1ex}
&S^1_3=\frac{r_1}{d_3},&&S^2_3=\frac{r_2}{d_3},&&S^3_3=\frac{r_3}{d_3}.
\hskip -2em
\endxalignat
$$
From \thetag{2.4} for entire numbers $p_1$, $p_2$, $p_3$, and $d_1$ in
\thetag{2.7} we derive the relationship
$$
(p_1)^2+(p_2)^2+ (p_3)^2=(d_1)^2.\hskip -2em
\tag2.8
$$
If four entire numbers satisfy the relationship \thetag{2.8}, we say that
they form {\bf Pytha\-gorean tetrad}. Each column in orthogonal matrix with
rational components is related with some Pythagorean tetrad of entire
numbers. Thus, in \thetag{2.7} we have three Pythagorean tetrads determined
by transition matrix $S$:
$$
\xalignat 3
&(p_1,p_2,p_3,d_1),&&(q_1,q_2,q_3,d_2),&&(r_1,r_2,r_3,d_3).\hskip -2em
\tag2.9
\endxalignat
$$
Pythagorean tetrads of entire numbers \thetag{2.9} are orthogonal to
each other in the sense of the following relationships:
$$
\align
&p_1\,q_1+p_2\,q_2+p_3\,q_3=0,\\
&p_1\,r_1+p_2\,r_2+p_3\,r_3=0,\\
&r_1\,q_1+r_2\,q_2+r_3\,q_3=0.
\endalign
$$
In order to determine an orthogonal matrix it's sufficient to have
two orthogonal Pythagorean tetrads, for instance, $(p_1,p_2,p_3,d_1)$
and $(q_1,q_2,q_3,d_2)$. Third Pythagorean tetrad then will be
determined by the relationships
$$
\xalignat 3
&\frac{r_1}{d_3}=-\frac{\vmatrix p_2 & p_3\\
\vspace{1ex}q_2 & q_3\endvmatrix
\vphantom{\vrule depth 12pt}}
{d_1\,d_2},
&&\frac{r_2}{d_3}=\frac{\vmatrix p_1 & p_3\\
\vspace{1ex}q_1 & q_3\endvmatrix
\vphantom{\vrule depth 12pt}}
{d_1\,d_2},
&&\frac{r_3}{d_3}=-\frac{\vmatrix p_1 & p_2\\
\vspace{1ex}q_1 & q_2\endvmatrix
\vphantom{\vrule depth 12pt}}
{d_1\,d_2}.
\endxalignat
$$
This is the consequence of the fact that third vector in orthonormal
bases (ONB's) are determined by vector product of first two vectors:
$$
\xalignat 2
&\bold e_3=-[\bold e_1,\,\bold e_2],
&&\bold h_3=[\bold h_1,\,\bold h_2].
\endxalignat
$$
The difference in sign here is due to the condition \thetag{2.6}, which
expresses difference in orientations of bases $\bold e_1$, $\bold e_2$,
$\bold e_3$ and $\bold h_1$, $\bold h_2$, $\bold h_3$.\par
     Returning to the problem on pyramid, we use the fact that vectors
$\overrightarrow{AF}$, $\overrightarrow{FH}$ and $\overrightarrow{HC}$
are collinear to base vectors $\bold e_1$, $\bold e_2$, $\bold e_3$:
$$
\xalignat 3
&\overrightarrow{AF}=\alpha_1\cdot\bold e_1,
&&\overrightarrow{FH}=\alpha_2\cdot\bold e_2,
&&\overrightarrow{HC}=\alpha_3\cdot\bold e_3.
\hskip -3em\tag2.10
\endxalignat
$$
Similarly, vectors $\overrightarrow{BG}$, $\overrightarrow{GK}$, and
$\overrightarrow{KC}$ are collinear to base vectors $\bold h_1$, $\bold
h_2$, $\bold h_3$:
$$
\xalignat 3
&\overrightarrow{BG}=\beta_1\cdot\bold h_1,
&&\overrightarrow{GK}=\beta_2\cdot\bold h_2,
&&\overrightarrow{KC}=\beta_3\cdot\bold h_3.
\hskip -3em
\tag2.11
\endxalignat
$$
Vector $\overrightarrow{BC}$ is collinear to vector $\bold e_3$, while
vector $\overrightarrow{AC}$ is collinear to vector $\bold h_3$:
$$
\xalignat 2
&\overrightarrow{BC}=\omega\cdot\bold e_3,
&&\overrightarrow{AC}=\sigma\cdot\bold h_3.
\hskip -3em
\tag2.12
\endxalignat
$$
Let's choose parameters $\omega$ and $\sigma$ in \thetag{2.12} to be
rational numbers, and then let's apply the relationships \thetag{2.1}.
As a result for vector $\overrightarrow{AC}$ we obtain two expansions:
$$
\aligned
\overrightarrow{AC}&=\overrightarrow{AF}+\overrightarrow{FH}+
\overrightarrow{HC},\\
\overrightarrow{AC}&=\sigma\cdot(S^1_3\,\bold e_1+S^2_3\,\bold e_2+
S^3_3\,\bold e_3).
\endaligned\hskip -3em
\tag2.13
$$
Substituting \thetag{2.10} into \thetag{2.13} and comparing two
expansions \thetag{2.13}, we obtain
$$
\xalignat 3
&\alpha_1=\sigma\,S^1_3,&&\alpha_2=\sigma\,S^2_3,&&\alpha_3=\sigma\,S^3_3.
\hskip -3em
\tag2.14
\endxalignat
$$
In a similar way, from \thetag{2.10} and \thetag{2.13} due to \thetag{2.1}
and due to the expansion $\overrightarrow{BC}=\overrightarrow{BG}
+\overrightarrow{GK}+\overrightarrow{KC}$ we can derive the following
three relationships:
$$
\xalignat 3
&\beta_1=\omega\,T^1_3,&&\beta_2=\omega\,T^2_3,&&\beta_3=\omega\,T^3_3.
\hskip -3em
\tag2.15
\endxalignat
$$
If components of transition matrix $S$ are rational numbers, then
components of inverse transition matrix $T=S^t$ are also rational.
Therefore from \thetag{2.14} and \thetag{2.15} we obtain rationality
of numeric coefficients $\alpha_2$, $\alpha_2$, $\alpha_3$ and $\beta_1$,
$\beta_2$, $\beta_3$ in \thetag{2.10} and \thetag{2.11}. This, in turn,
provides rationality of lengths of segments $[AC]$, $[BC]$, $[AK]$,
$[CK]$, $[BH]$, $[CH]$, $[BG]$, $[GK]$, $[AF]$, and $[FH]$.\par
     Let's consider the vector $\overrightarrow{FG}$, length of which
is a final result in the problem on pyramid. For this vector we have
an expansion:
$$
\overrightarrow{FG}=\overrightarrow{FH}+\overrightarrow{HC}-
\overrightarrow{KC}-\overrightarrow{GK}.\hskip -3em
\tag2.16
$$
Let's substitute \thetag{2.10} and \thetag{2.11} into the expansion
\thetag{2.16}. This yields
$$
\align
\overrightarrow{FG}&=\alpha_2\,\bold e_2+\alpha_3\,\bold e_3-\beta_2\,
\bold h_2-\beta_3\,\bold h_3=\alpha_2\,\bold e_2+\alpha_3\,\bold e_3\,-\\
&-\,\beta_2\,(S^1_2\,\bold e_1+S^2_2\,\bold e_2+S^3_2\,\bold e_3)
-\beta_3\,(S^1_3\,\bold e_1+S^2_3\,\bold e_2+S^3_3\,\bold e_3).
\endalign
$$
Now it's easy to see that vector $\overrightarrow{FG}$ has rational
coordinates in orthonormal base (ONB) formed by vectors $\bold
e_1$, $\bold e_2$, $\bold e_3$. Therefore its length, in the worst
case, is simplest irrational number obtained as a {\bf square root of
rational number}. The same is true for lengths of segments $[AB]$,
$[AM]$, $[BM]$, $[CM]$, $[KL]$, $[HL]$, $[AH]$, $[BK]$, as well as
for the lengths of segments $[F\tilde F]$, $[H\tilde F]$, $[G\tilde G]$,
and $[K\tilde G]$. For $|AB|$ this follows from the equality
$\overrightarrow{AB}=\overrightarrow{AC}-\overrightarrow{BC}$. Further
we use formulas \thetag{1.2}, \thetag{1.3}, \thetag{1.5}, \thetag{1.4};
then formulas \thetag{1.6} and \thetag{1.7}. Main conclusion that we
draw from what was said above is the following: {\bf orthogonal matrices
with rational components} give the {\bf algorithm} for choosing
numeric values of parameters $a$, $b$, $c$, $f$, $g$ in the problem on
pyramid such that we get simple final result in this problem and simple
results in all intermediate calculations.
\head
3. Constructing orthogonal matrices\\with rational components.
\endhead
     Constructing orthogonal matrices with rational components is a
separate problem. First we consider regular algorithm for constructing
such matrices. It is based on elementary rotations. Let's consider
three entire numbers $p_1$, $p_2$, $d$, and suppose that they are bound
by the relationship
$$
(p_1)^2+(p_2)^2=d^2.\hskip -2em
\tag3.1
$$
Such numbers form {\bf Pythagorean triad}. In contrast to Pythagorean
tetrads they are well-known. There is a regular algorithm for constructing
all Pythagorean triads (see, for instance, \cite{1}). If $\tau$ is a
greatest common divisor of $p_1$, $p_2$, and $q$, then $p_1=\tau\cdot\tilde
p_1$, \ $p_2=\tau\cdot\tilde p_2$, \ $d=\tau\cdot\tilde d$. From \thetag{3.1}
we derive
$$
(\tilde p_1)^2+(\tilde p_2)^2=\tilde d^2.
$$
If $\tilde p_1$ is even and $\tilde p_2$ is odd, then $\tilde d$ is odd.
According to the regular algorithm described in \cite{1}, in this case we
have the following expressions:
$$
\align
&\tilde p_1=2\,(m^2+m-n^2-n),\\
&\tilde p_2=4\,m\,n+2\,m+2\,n+1,\\
&\tilde d=2\,(m^2+m+n^2+n)+1.
\endalign
$$
Here $m$ and $n$ are two arbitrary entire numbers. So Pythagorean triads
are parameterized by three arbitrary entire numbers: $m$, $n$, and $\tau$.
     Suppose that we have some nonzero Pythagorean triad $(p_1,p_2,d)$.
Then we can consider two rational numbers $p_1/d$ and $p_2/d$, sum of
their squares being equal to unity:
$$
\left(\frac{p_1}{d}\right)^2+\left(\frac{p_2}{d}\right)^2=1.
$$
Hence in half-open interval $[0,\,2\pi)$ there exists some angle $\varphi$
such that
$$
\xalignat 2
&\cos\varphi=\frac{p_1}{d},&&\sin\varphi=\frac{p_2}{d}.\hskip -2em
\tag3.2
\endxalignat
$$
Angle $\varphi$ in \thetag{3.2} is uniquely determined by numbers
$p_1$, $p_2$, and $d$ forming Pythago\-rean triad $(p_1,p_2,q)$.
Let's use this angle in order to define four matrices:
$$
\xalignat 2
&S^\sqbi{x}_\varphi=\Vmatrix
\cos\varphi & \sin\varphi & 0\\
\vspace{2ex}
-\sin\varphi & \cos\varphi & 0\\
\vspace{2ex}
0 & 0 & 1 \endVmatrix,
&&S^\sqbi{y}_\varphi=\Vmatrix
\cos\varphi & 0 & \sin\varphi \\
\vspace{2ex}
0 & 1 & 0\\
\vspace{2ex}
-\sin\varphi & 0 & \cos\varphi
\endVmatrix,\hskip -3em\\
&&&\tag3.3\\
\hskip -3em
&S^\sqbi{z}_\varphi=\Vmatrix
1 & 0 & 0\\
\vspace{2ex}
0 & \cos\varphi & \sin\varphi\\
\vspace{2ex}
0 &-\sin\varphi & \cos\varphi
\endVmatrix,
&&S^*=\Vmatrix
-1 & 0 & 0 \\
\vspace{2ex}
0 & -1 & 0\\
\vspace{2ex}
0 & 0 & -1\endVmatrix.\hskip -3em
\endxalignat
$$
Matrices $S^\sqbi{x}_\varphi$, $S^\sqbi{y}_\varphi$, $S^\sqbi{z}_\varphi$ 
are geometrically interpreted as matrices of elementary rotations to the
angle $\varphi$ around coordinate axes. They arise as transition matrices
in \thetag{2.1} in that case when base $\bold h_1$, $\bold h_2$, $\bold h_3$
is got from base $\bold e_1$, $\bold e_2$, $\bold e_3$ by one of such
elementary rotations. Matrix $S^*$ in \thetag{3.3} is interpreted as a
matrix of inversion. It arises as transition matrix in \thetag{2.1} when
$$
\xalignat 3
&\bold h_1=-\bold e_1,&&\bold h_2=-\bold e_2,&&\bold h_3=-\bold e_3.
\endxalignat
$$\par
     All four matrices \thetag{3.3} are orthogonal. This can be checked
by substituting them into \thetag{2.3}. If $\varphi$ is determined by the
relationships \thetag{3.2}, then all components of matrices \thetag{3.3}
are rational. Product of two orthogonal matrices is an orthogonal matrix
(it's well-known that such matrices form a group). Therefore, choosing
$n$ Pythagorean triads and determining angles $\varphi_1,\,\ldots,\,
\varphi_n$ by relationships \thetag{3.2}, we can consider the following
product of corresponding matrices \thetag{3.3}:
$$
S=(S^*)^\varepsilon\cdot\prod^n_{i=1}
(S^\sqbi{x}_{\varphi_i})^{\alpha_i}
\cdot(S^\sqbi{y}_{\varphi_i})^{\beta_i}
\cdot(S^\sqbi{z}_{\varphi_i})^{\gamma_i}.
\hskip -2em
\tag3.4
$$
Here $\varepsilon$, $\alpha_i$, $\beta_i$, $\gamma_i$ are entire numbers
either equal to zero or to unity. Formula \thetag{3.4} gives an algorithm
for constructing orthogonal matrices with rational components in dimension
$3$. We shall call it a {\bf regular algorithm}.
\head
4. Some generalizations and open questions.
\endhead
    Leonard Euler considered a class of matrices, which is a little more
wide than class determined by the relationships \thetag{2.3}. Euler's class
is formed by matrices $S$ with entire components that satisfy the following
condition:
$$
S^t\cdot S=N\cdot E.\hskip -2em
\tag4.1
$$
Here $N$ is some positive entire number. Matrices of Euler's class are
called {\bf entire orthogonal matrices}, number $N$ is called a {\bf norm
of orthogonality}. If $N$ is a square of entire number, i\.~e\. if $N=M^2$,
then matrix $M^{-1}\cdot S$ is an orthogonal matrix with rational
components in the sense of standard definition by formulas \thetag{2.3}.
Leonard Euler has suggested an algorithm for constructing entire
orthogonal matrices in the dimensions $3$ and $4$. His algorithm is
described in book \cite{2}. In papers \cite{3--5} Euler's algorithm was
generalized for $n\times n$ matrices in arbitrary dimension $n$. Due to
the existence of two algorithms we have a series of quite natural
questions.
\roster
\rosteritemwd=2pt
\item"---" How do Euler's algorithm relate with regular algorithm, which
      is expressed by above formula \thetag{3.4}\,?
\item"---" Can we construct an arbitrary orthogonal matrix with rational
      components by Euler's algorithm\,?
\item"---" Is there the expansion \thetag{3.4} for an arbitrary orthogonal
      matrix with rational components, i\.~e\. can it be constructed by
      regular algorithm\,?
\endroster
Answers to these questions are unknown to the author of this paper. Author
will be grateful for any information concerning subject of this paper.
\head
5. Acknowledgements.
\endhead
     Paper was reported at the conference of Soros Educators in
1998 (Beloretsk, Russia). Author is grateful to Dr.~V.~G.~Khazankin
for invitation to this conference. Author is also grateful to George
Soros foundation (Open Society Institute) for financial support in
1998 (grant No\.~d98-943).

\enddocument
\end